
\input amssym.def
\input amssym
\magnification=1100
\baselineskip = 0.25truein
\lineskiplimit = 0.01truein
\lineskip = 0.01truein
\vsize = 8.5truein
\voffset = 0.2truein
\parskip = 0.10truein
\parindent = 0.3truein
\settabs 12 \columns
\hsize = 5.4truein
\hoffset = 0.4truein

\setbox\strutbox=\hbox{%
\vrule height .708\baselineskip
depth .292\baselineskip
width 0pt}
\font\caps=cmcsc10

\font\bigtenrm=cmr10 at 14pt

\def\sqr#1#2{{\vcenter{\vbox{\hrule height.#2pt
\hbox{\vrule width.#2pt height#1pt \kern#1pt
\vrule width.#2pt}
\hrule height.#2pt}}}}
\def\square{\mathchoice\sqr46\sqr46\sqr{3.1}6\sqr{2.3}4}

\centerline{\bigtenrm THE ASYMPTOTIC BEHAVIOUR}
\centerline{\bigtenrm OF HEEGAARD GENUS}

\tenrm
\vskip 14pt
\centerline{MARC LACKENBY}
\vskip 18pt

\tenrm
\centerline{\caps 1. Introduction}
\vskip 6pt

Heegaard splittings have recently been shown to be related
to a number of important conjectures in 3-manifold theory:
the virtually Haken conjecture, the positive virtual
$b_1$ conjecture and the virtually fibred conjecture [3].
Of particular importance is the rate at which the Heegaard
genus of finite-sheeted covering spaces grows as a function of their
degree. This was encoded in the following definitions.

Let $M$ be a compact orientable 3-manifold. Let $\chi_-^h(M)$
be the negative of the maximal Euler characteristic of a
Heegaard surface for $M$. Let $\chi_-^{sh}(M)$ be
the negative of the maximal Euler characteristic of
a strongly irreducible Heegaard surface for $M$,
or infinity if such a surface does not exist.
Define the {\sl infimal Heegaard gradient} of $M$
to be
$$\inf_i \{ {\chi_-^h(M_i) \over d_i} : M_i 
\hbox{ is a degree } d_i \hbox{ cover of } M \}
.$$
The {\sl infimal strong Heegaard gradient} of $M$ is
$$\liminf_i \{ {\chi_-^{sh}(M_i) \over d_i} : M_i 
\hbox{ is a degree } d_i \hbox{ cover of } M \}.$$

The following conjectures were put forward in [3].
According to Theorem 1.7 of [3], either of these conjectures,
together with a conjecture of Lubotzky and Sarnak [4]
about the failure of Property ($\tau$) for
hyperbolic 3-manifolds, would imply the
virtually Haken conjecture for hyperbolic 3-manifolds.

\noindent {\bf Heegaard gradient conjecture.} {\sl A compact orientable 
hyperbolic 3-manifold has zero infimal Heegaard gradient if and only if it
virtually fibres over the circle.}

\noindent {\bf Strong Heegaard gradient conjecture.} 
{\sl Any closed orientable hyperbolic
3-manifold has positive infimal strong Heegaard gradient.}

Some evidence for these conjectures was presented
in [3]. More precisely, suitably phrased versions of
these conjectures were shown to be true
when one restricts attention to cyclic covers
dual to a non-trivial element of $H_2(M, \partial M)$, to reducible
manifolds and (in the case of the strong Heegaard
gradient conjecture) to congruence covers of arithmetic hyperbolic
3-manifolds. 

A less quantitative version of the conjectures
is simply that $\chi_-^{sh}(M_i)$ cannot grow
too slowly as a function of $d_i$, and that
if $\chi_-^{h}(M_i)$ does grow sufficiently slowly,
then $M$ is virtually fibred. These expectations are 
confirmed in the following result, which is
the main theorem of this paper.

\noindent {\bf Theorem 1.} {\sl Let $M$ be a closed
orientable 3-manifold that admits a negatively
curved Riemannian metric. Let $\{ M_i \rightarrow M \}$
be a collection of finite regular covers with degree $d_i$.
\item{(1)} If $\chi_-^h(M_i) / \sqrt{d_i} \rightarrow 0$,
then $b_1(M_i) > 0$ for all sufficiently large $i$.
\item{(2)} $\chi_-^{sh}(M_i) / \sqrt{d_i}$ is bounded away from zero.
\item{(3)} If $\chi_-^h(M_i) / \root 4 \of {d_i} \rightarrow 0$,
then $M_i$ fibres over the circle for all sufficiently large $i$.}

A slightly weaker form of Theorem 1(1)
appeared in [3] as Corollary 1.4, with essentially the same proof.
It is included here in order to emphasise its connection
to the other two results. 

The following corollary of Theorem 1(3) gives a necessary
and sufficient condition for $M$ to be virtually fibred
in terms of the Heegaard genus of its finite covers.
We say that a collection $\{ M_i \rightarrow M \}$ of
finite covers has {\sl bounded irregularity} if
the normalisers of $\pi_1 M_i$ in $\pi_1 M$ have
bounded index in $\pi_1 M$.

\noindent {\bf Corollary 2.} {\sl Let $M$ be a closed orientable
3-manifold with a negatively curved Riemannian metric, and
let $\{ M_i \rightarrow M \}$ be its finite-sheeted
covers with degree $d_i$. Then the following are equivalent:
\item{(1)} $M_i$ is fibred for infinitely many $i$;
\item{(2)} in some subsequence with bounded
irregularity, $\chi_-^h(M_i)$ is bounded;
\item{(3)} in some subsequence with bounded
irregularity, $\chi_-^h(M_i) / \root 4 \of {d_i}
\rightarrow 0$.

}

\noindent {\sl Proof.} $(1) \Rightarrow (2)$. If some
finite-sheeted cover $\tilde M$ is fibred, then so is any
finite cyclic cover $M_i$ of $\tilde M$ dual to the fibre. The
normaliser $N(\pi_1M_i)$ of $\pi_1M_i$ in $\pi_1 M$
contains $\pi_1\tilde M$,
so $[\pi_1M : N(\pi_1M_i)]$ is bounded and hence
these covers have bounded irregularity. Also, 
$\chi_-^h(M_i)$ is bounded by twice the modulus of the Euler
characteristic of the fibre, plus four.

$(2) \Rightarrow (3)$. This is trivial, since $d_i$ must tend to infinity.

$(3) \Rightarrow (1)$. Since $[\pi_1M : N(\pi_1M_i)]$ is
bounded in this subcollection, we may pass to a further subsequence where
$N(\pi_1M_i)$ is a fixed subgroup of $\pi_1M$. Let $\tilde M$ be the
finite-sheeted cover of $M$ corresponding to this
subgroup. Then, the covers $\{ M_i \rightarrow M \}$
in this subsequence give a collection
$\{ M_i \rightarrow \tilde M \}$ of finite regular covers such that
$\chi_-^h(M_i) / [\pi_1 \tilde M : \pi_1 M_i ]^{1/4}
\rightarrow 0$. By Theorem 1(3), $M_i$
is fibred for all sufficiently large $i$. $\square$

\vskip 18pt
\centerline {\caps 2. Background material}
\vskip 6pt

\noindent {\caps Generalised Heegaard splittings}
\vskip 6pt

A Heegaard splitting of a closed orientable 3-manifold
can be viewed as arising from a handle structure. If one builds the
manifold by starting with a single 0-handle, then
attaching some 1-handles, then some 2-handles
and then a 3-handle, the manifold obtained
after attaching the 0- and 1-handles is handlebody,
as is the closure of its complement. Thus, the boundary
of this submanifold is a Heegaard surface.
Generalised Heegaard splittings arise from more
general handle structures: one starts with some
0-handles, then adds some 1-handles, then some 2-handles,
then 1-handles, and so on, in an alternating
fashion, ending with some 3-handles. One then
considers the manifold embedded in $M$ consisting
of the 0-handles and the first $j$ batches of
1- and 2-handles. Let $F_j$ be the boundary of
this manifold, but discarding any 2-sphere
components that bound 0- or 3-handles. After
a small isotopy, so that these surfaces are all
disjoint, they divide $M$ into compression bodies.
In fact, the surfaces $\{ F_j : j \hbox{ odd} \}$
form Heegaard surfaces for the manifold
$M - \bigcup \{ F_j : j \hbox{ even} \}$.
We term the surfaces $F_j$ {\sl even} or {\sl odd},
depending on the parity of $j$.

More details about generalised Heegaard splittings
can be found in [7] and [6]. The following theorem
summarises some of the results from [7].

\noindent {\bf Theorem 3.} {\sl From any
minimal genus Heegaard surface $F$ for a closed orientable
irreducible 3-manifold $M$, other than $S^3$, one can construct
a generalised Heegaard splitting $\{ F_1, \dots, F_n \}$
in $M$ with the following properties:
\item{1.} $F_j$ is incompressible and has no 2-sphere
components, for each even $j$;
\item{2.} $F_j$ is strongly irreducible for each odd $j$;
\item{3.} $F_j$ and $F_{j+1}$ are not parallel for
any $j$;
\item{4.} $|\chi(F_j)| \leq |\chi(F)|$ for each $j$;
\item{5.} $|\chi(F)| = \sum (-1)^j \chi(F_j)$.

}

\noindent {\bf Corollary 4.} {\sl Let $M$, $F$ and
$\{ F_1, \dots, F_n \}$ be as in Theorem 3. Suppose
that, in addition, $M$ is not a lens space. Let $\overline F$
be the surface obtained from $\bigcup_j F_j$ by
replacing any components that are parallel by a single
component. Then
\item{1.} $|\chi(\overline F)| \leq |\chi(\bigcup_j F_j)|
< |\chi(F)|^2$;
\item{2.} $\overline F$ has at most ${3 \over 2} |\chi(F)|$ components.

}

\noindent {\sl Proof.} Note first that no component of
$\bigcup_j F_j$ is a 2-sphere. When $j$ is even, this is (1)
of Theorem 3. The same is true when $j$ is odd, since the
odd surfaces form Heegaard surfaces for the complement
of the even surfaces, and $M$ is not $S^3$.
Hence, none of the compression bodies $H$
in the complement of $\bigcup_j F_j$ is a 3-ball.

We claim also that no $H$ is a solid
torus. For if it were, consider the compression
body to which it is adjacent. If this were a product,
some even surface would be compressible, contradicting
(1) of Theorem 3. However, if it was not a product,
then it is a solid torus, possibly with punctures. If
it has no punctures, then $M$ is a lens space, contrary to
assumption. If it
does, then some component of an even surface would be
a 2-sphere, again contradicting (1). 
This proves the claim.

We expand (5) of Theorem 3 as follows:
$$ |\chi(F)| = {-\chi(F_1) \over 2} + {\chi(F_2) - \chi(F_1) \over 2} + 
{\chi(F_2) - \chi(F_3) \over 2}
+ \dots + {-\chi(F_n) \over 2}. \eqno{(\ast)}$$
For any compression body
$H$, other than a 3-ball, with negative boundary $\partial_- H$ and
positive boundary $\partial_+ H$, 
$\chi(\partial_- H) - \chi(\partial_+ H)$
is even and non-negative. It is zero
if and only if $H$ is a product or a solid torus.
Since $F_j$ and $F_{j+1}$ are not parallel
for any $j$, each term in ($\ast$) is
therefore at least one. So, $n+1$, the
number of terms on the right-hand side of ($\ast$), is at most
$|\chi(F)|$. Hence,
$$|\chi(\bigcup_j F_j)| = \sum_j |\chi(F_j)|
\leq n |\chi(F)| < |\chi(F)|^2.$$
The inequality $|\chi(\overline F)| \leq |\chi(\bigcup_j F_j)|$
simply follows from the fact that we discard some
components of $\bigcup_j F_j$ to form $\overline F$.
This proves (1).

Now it is trivial to check that, for any compression
body $H$, other than a 3-ball, solid torus
or product, $|\partial H| \leq {3 \over 2}(\chi(\partial_-H) 
- \chi(\partial_+ H))$. The number of components
of $\overline F$ is half the sum, over all complementary
regions $H$ of $\bigcup_j F_j$ that are not products,
of $|\partial H|$. This is at most ${3 \over 4}(\chi(\partial_-H) 
- \chi(\partial_+ H))$. But the sum, over all complementary
regions $H$ of $\bigcup_j F_j$, of
${1 \over 2}(\chi(\partial_-H) - \chi(\partial_+ H))$ is
the right-hand side of ($\ast$). Thus, we deduce that the
number of components of $\overline F$ is at most
${3 \over 2} |\chi(F)|$, proving (2). $\square$

\vskip 6pt
\noindent {\caps Realisation as minimal surfaces}
\vskip 6pt

One advantage of using generalised Heegaard splittings
satisfying (1) and (2) of Theorem 3 is that minimal
surfaces then play a r\^ole in the theory. The following
theorem of Freedman, Hass and Scott [2] applies to the
even surfaces.

\noindent {\bf Theorem 5.} {\sl Let $S$ be an orientable
embedded incompressible surface in a closed orientable irreducible
Riemannian 3-manifold. Suppose that no two components
of $S$ are parallel, and that no component is a 2-sphere. Then there is an ambient isotopy
of $S$ so that afterwards each component is either a
least area, minimal surface or the boundary of a regular
neighbourhood of an embedded, least area, minimal non-orientable
surface.}

We will apply the above result to the incompressible components
of $\overline F$. If we cut $M$ along these components,
the remaining components form strongly
irreducible Heegaard surfaces for the complementary
regions. A theorem of Pitts and Rubinstein [5] now
applies.

\noindent {\bf Theorem 6.} {\sl Let $S_1$ be a (possibly empty) 
embedded stable minimal surface in a closed orientable irreducible
Riemannian 3-manifold $M$ with a bumpy metric. Let $S_2$
be a strongly irreducible Heegaard surface for a
complementary region of $S_1$. Then there is an ambient
isotopy, leaving $S_1$ fixed, taking $S_2$ to a
minimal surface, or to the boundary of a regular
neighbourhood of a minimal embedded non-orientable surface, with
a tube attached that is vertical in the $I$-bundle 
structure on this neighbourhood.}

Bumpy metrics were defined by White in [8]. After
a small perturbation, any Riemannian metric can be made bumpy.
Then we may ambient isotope $\overline F$ so that
each component is as described in Theorems 5 and 6.

We will need some parts of the proof of Theorem 6,
and not just its statement. Let $X$ be the component
of $M - S_1$ containing $S_2$. Then, as a Heegaard
surface, $S_2$ determines a sweepout of $X$. In any
sweepout, there is a surface of maximum area, although
it need not be unique. Let $a$ be the infimum,
over all sweepouts in this equivalence class, of this
maximum area. Then Pitts and Rubinstein showed
that there is a sequence of sweepouts, whose
maximal area surfaces tend to an embedded minimal surface,
and that the area of these surfaces tends to $a$.
This minimal surface, or its orientable double
cover if it is non-orientable, is isotopic to $S_2$
or to a surface obtained by compressing $S_2$.

Now, when $M$ is negatively curved, one may use
Gauss-Bonnet to bound the area of this surface.
Suppose that $\kappa < 0$ is the supremum of the
sectional curvatures of $M$. Then, as the
surface is minimal, its sectional curvature is
at most $\kappa$. Hence, by Gauss-Bonnet,
its area is at most $2 \pi |\chi(S_2)| /
|\kappa|$. Thus, we have the following result.

\noindent {\bf Addendum 7.} {\sl Let $S_1$,
$S_2$ and $M$ be as in Theorem 6. Suppose that
the sectional curvature of $M$ is at most
$\kappa < 0$. Then, for each $\epsilon > 0$,
there is a sweepout of the component of
$M - S_1$ containing $S_2$, equivalent
to the sweepout determined by $S_2$, so that each surface
in this sweepout has area at most
$(2 \pi |\chi(S_2)|/ |\kappa|) + \epsilon$.}

One has a good deal of geometric control over
minimal surfaces when $M$ is negatively curved.
As observed above, their area is bounded
in terms of their Euler characteristic and
the supremal sectional curvature of $M$. In fact,
by ruling out the existence of long thin
tubes in the surface, one has the following.

\noindent {\bf Theorem 8.} {\sl There is
a function $f \colon {\Bbb R} \times {\Bbb R}
\rightarrow {\Bbb R}$ with the following
property. Let $M$ be a Riemannian 3-manifold,
whose injectivity radius is at least $\epsilon /2 > 0$,
and whose sectional curvature is at most
$\kappa < 0$. Let $S$ be a closed minimal
surface in $M$. Then there is a collection of
at most $f(\kappa, \epsilon) |\chi(S)|$
points in $S$, such that the balls of radius
$f(\kappa, \epsilon)$ about these points cover $S$. 
(Here, distance is
measured using the path metric on $S$.)}

This is proved in Proposition 6.1 of [3]. More
precisely, formulas (1) and (2) there give the
result.

\vskip 6pt
\noindent {\caps The Cheeger constant of manifolds and graphs}
\vskip 6pt

The {\sl Cheeger constant} of a compact Riemannian
manifold $M$ is defined to be
$$h(M) = \inf \left\{ { {\rm Area}(S) \over \min \{ 
{\rm Volume}(M_1), {\rm Volume}(M_2) \} } \right\},$$
where $S$ ranges over all embedded codimension one
submanifolds that divide $M$ into $M_1$ and $M_2$.

A central theme of [3] is that the Cheeger constant
of a 3-manifold and its Heegaard splittings are
intimately related. One example of this phenomenon is the
following result.

\noindent {\bf Theorem 9.} {\sl Let $M$ be a closed
Riemannian 3-manifold. Let $\kappa <0$
be the supremum of its sectional curvatures. Then
$$h(M) \leq { 4 \pi \, \chi_-^h(M) \over |\kappa|
{\rm Volume}(M)}.$$

}

This is essentially Theorem 4.1 of [3]. However, there,
$\chi_-^h(M)$ is replaced by $c_+(M)$, which is an
invariant defined in terms of the generalised Heegaard
splittings of $M$. But the above inequality follows
from an identical argument. We briefly summarise the 
proof.

From a minimal genus Heegaard splitting of $M$, construct
a generalised Heegaard splitting $\{ F_1, \dots, F_n \}$
satisfying (1) to (5) of Theorem 3. Let
$\overline F$ be the surface obtained from $\bigcup_j F_j$
by discarding multiple copies of parallel components.
Apply the isotopy of Theorem 5 to the incompressible components
of $\overline F$. Each complementary
region corresponds to a component of the complement of
the even surfaces, and therefore contains a component of
some odd surface $F_j$. Label this region with the integer
$j$, and let $M_j$ be the union of the regions labelled $j$.
There is some odd $j$ such that the volumes of
$M_1 \cup \dots \cup M_{j-2}$ and $M_{j+2} \cup \dots \cup M_n$
are each at most half the volume of $M$. Now, $F_j \cap M_j$
forms a strongly irreducible Heegaard surface for $M_j$. Applying Addendum 7,
we find for each $\epsilon > 0$, a sweepout of $M_j$,
equivalent to that determined by $F_j$, by surfaces with area
at most $(2 \pi |\chi(F_j)| / |\kappa|) + \epsilon$. But,
$|\chi(F_j)| \leq \chi_-^h(M)$, by (4) of Theorem 3.
Some surface in this sweepout divides $M$ into two
parts of equal volume. So, as $\epsilon$ was arbitrary,
$$h(M) \leq { 4 \pi \, \chi_-^h(M) \over |\kappa|
{\rm Volume}(M)}.$$

In this paper, we will consider the Cheeger constants of
regular finite-sheeted covering spaces $M_i$ of $M$. Here,
$M_i$ is given the Riemannian metric lifted from $M$. It
is possible to estimate $h(M_i)$ in terms of graph-theoretic
data, as follows.

By analogy with the Cheeger constant for a Riemannian
manifold, one can define the {\sl Cheeger constant} $h(X)$ of a finite
graph $X$. If $A$ is a subset of the vertex set $V(X)$,
$\partial A$ denotes those edges with precisely one endpoint
in $A$. Then $h(X)$ is defined to be
$$\inf \left \{ {|\partial A| \over |A|}: A \subset V(X)
\hbox{ and } 0 < |A| \leq |V(X)|/2 \right \}.$$

\noindent {\bf Proposition 10.} {\sl Let $M$ be
a compact Riemannian manifold. Let ${\cal X}$
be a finite set of generators for $\pi_1M$. Then there
is a constant $k_1 \geq 1$ with the following property.
If $X_i$ is the Cayley graph of $\pi_1 M / \pi_1 M_i$
with respect to the generators ${\cal X}$, then
$$k_1^{-1} \, h(X_i) \leq h(M_i) \leq k_1 \, h(X_i).$$}

This is essentially contained in [1], but we outline a proof. 
Lemma 2.3 of [3] states that, if ${\cal X}$ and ${\cal X}'$ are two finite
sets of generators for $\pi_1 M$, then there is a constant
$k \geq 1$ with the following property. If $X_i$ and
$X'_i$ are the Cayley graphs of $\pi_1 M / \pi_1 M_i$
with respect to ${\cal X}$ and ${\cal X}'$, then
$$k^{-1} \, h(X_i) \leq h(X'_i) \leq k \, h(X_i).$$
Thus, for the purposes of proving Proposition 10, we
are free to choose ${\cal X}$. We do this as follows.
We pick a connected fundamental domain in the universal
cover of $M$. The translates of this domain to which it
is adjacent correspond to a finite set ${\cal X}$
of generators for $\pi_1 M$. There is an induced
fundamental domain in any finite regular cover $M_i$
of $M$. Its translates are in one-one correspondence
with the group $\pi_1 M / \pi_1 M_i$. Two translates
are adjacent if and only if one is obtained from the
other by right-multiplication by an element in ${\cal X}$.
Thus, the Cayley graph $X_i$ should be viewed as a
coarse approximation to $M_i$. Any subset $A$ of
$V(X_i)$, as in the definition of $h(X_i)$, therefore
determines a decomposition of $M_i$. After a further
modification, we may assume that this is along a codimension
one submanifold. The existence
of a constant $k_1$ such that $h(M_i) \leq k_1 \, h(X_i)$
is then clear. The other inequality is more difficult to
establish. One needs to control the geometry of
a codimension one submanifold $S$ in $M_i$ that is
arbitrarily close to realising the Cheeger constant
of $M_i$. This is achieved in the proof of Lemma 2 of [1].

\vskip 6pt
\noindent {\caps Constructing non-trivial cocycles}
\vskip 6pt

Some new machinery has been developed in [3] that
gives necessary and sufficient conditions on a finitely
presented group to have finite index subgroups with
infinite abelianisation. We describe some of the ideas
behind this now.

Let $C$ be a finite cell complex with a single
0-cell and in which every 2-cell is a triangle.
Let $G$ be its fundamental group, and let
${\cal X}$ be the generators arising from the
1-cells. Associated with any finite index
normal subgroup $H_i$ of $G$, there is a 
finite-sheeted covering space $C_i$ of $C$.
Its 1-skeleton $X_i$ is the Cayley graph of
$G/H_i$ with respect to ${\cal X}$. The following
theorem is an expanded form of Lemma 2.4 of [3]
and has exactly the same proof.
It will play a key r\^ole in this paper.

\noindent {\bf Theorem 11.} {\sl Suppose that
$h(X_i) < \sqrt{ 2 / (3 |V(X_i)|)}$.
Let $A$ be any non-empty subset of
$V(X_i)$ such that $|\partial A|/|A| = h(X_i)$
and $|A| \leq |V(X_i)|/2$. Then there is a
1-cocycle $c$ on $C_i$ that is not a 
coboundary. Its support is a subset of the edges
of $\partial A$, and it takes values in
$\{ -1, 0, 1 \}$. As a consequence, $H_i$
has infinite abelianisation.}

\vskip 18pt
\centerline {\caps 3. The proof of the main theorem}
\vskip 6pt

We start with a closed orientable 3-manifold $M$ admitting
a negatively curved Riemannian metric. After a small perturbation, we may assume
that the metric is bumpy. Let $\kappa < 0$
be the supremum of its sectional curvatures.
Pick a 1-vertex triangulation $T$ of $M$. The
edges of $T$, when oriented in some way, form a set ${\cal X}$ of generators
for $\pi_1(M)$. Let $K$ be the 2-skeleton of
the complex dual to $T$. 

We will consider a collection
$\{ M_i \rightarrow M \}$ of finite regular covers
of $M$, having the properties of Theorem 1.
In particular, we will assume (at least) that
$\chi_-^h(M_i) / \sqrt d_i \rightarrow 0$.
(Note that this is justified when proving Theorem 1(2),
by passing to a subsequence, and using the fact that
$\chi_-^h(M_i) \leq \chi_-^{sh}(M_i)$.)
The triangulation $T$ and 2-complex $K$ lift to $T_i$
and $K_i$, say, in $M_i$. The 1-skeleton of $T_i$
forms the Cayley graph $X_i$ of $G_i = \pi_1 M / \pi_1 M_i$ with 
respect to ${\cal X}$.

According to Theorem 9,
$$h(M_i) \leq {4 \pi \over |\kappa|} {\chi_-^h(M_i) \over {\rm Volume}(M_i)} = 
{4 \pi \over |\kappa| {\rm Volume}(M)} {\chi_-^h(M_i) \over d_i},$$
By Proposition 10, there is a constant $k_1 \geq 1$ independent of $i$ such that
$h(X_i) \leq k_1 h(M_i)$. Setting
$$k_2 = {4 \pi k_1 \over |\kappa| {\rm Volume}(M)},$$
we deduce that
$$h(X_i) \leq k_2 {\chi_-^h(M_i) \over d_i}.$$ 

Let $V(X_i)$ be the vertex set of $X_i$. Let $A$ be a non-empty
subset of $V(X_i)$ such that $|\partial A|/|A|
= h(X_i)$ and $|A| \leq |V(X_i)| / 2 = d_i/2$.
By Theorem 11, when $h(X_i) < \sqrt{2/(3 d_i)}$,
$T_i$ admits a 1-cocycle $c$ that is not a coboundary.
Since $h(X_i) \leq k_2 \chi_-^h(M_i) /d_i$, and we are assuming
(at least) that $\chi_-^h(M_i) / \sqrt{d_i} \rightarrow 0$,
then such a cocycle exists for all sufficiently large $i$.
This establishes (1) of the Theorem 1.

Theorem 11 states that $c$ takes values in $\{-1,0,1\}$, and
its support is a subset of the edges of $\partial A$. 
Dual to this cocycle is a transversely oriented normal surface $S$ in $T_i$
which is homologically non-trivial. Remove any
2-sphere components from $S$. This is still homologically non-trivial,
since all 2-spheres in $M_i$ are inessential, as $M_i$
is negatively curved.
The intersection of $S$ with the 2-skeleton of $T_i$ is a graph in $S$ whose
complementary regions are triangles and squares.
Let $V(S)$ and $E(S)$ be its vertices and
edges. Its vertices are in one-one correspondence with
the edges of $T_i$ in the support of $c$. So, 
$$|V(S)| \leq |\partial A| = |A| h(X_i) \leq d_i h(X_i) / 2 
\leq k_2 \chi_-^h(M_i) /2.$$
The valence of each vertex is at most the maximal
valence of an edge in $T$, $k_3$, say.
So, $$|\chi(S)| < |E(S)| \leq |V(S)| k_3 /2 \leq 
k_2  k_3 \chi_-^h(M_i) /4.$$
Setting $k_4 = k_2 k_3 /4 $, we have deduced
the existence of a homologically non-trivial, transversely 
oriented, properly embedded surface $S$
with $|\chi(S)| \leq k_4 \chi_-^h(M_i)$ and with no
2-sphere components. By compressing $S$ and removing
components if necessary, we may assume that $S$
is also incompressible and connected. 
Thus, we have proved the following result.

\noindent {\bf Theorem 12.} {\sl Let $M$ be a closed
orientable 3-manifold with a negatively curved
Riemannian metric. Then there is a constant $k_4 > 0$
with the following property. Let $\{ M_i \rightarrow M \}$
be a collection of finite regular covers, with degree
$d_i$. If $\chi_-^h(M_i) / \sqrt{d_i} \rightarrow 0$, then,
for all sufficiently large $i$, $M_i$ contains an embedded,
connected, oriented, incompressible, homologically non-trivial
surface $S$ such that $|\chi(S)| \leq k_4 \chi_-^h(M_i)$.}

By a theorem of Freedman, Hass and Scott [2] (Theorem 5 in this paper), 
there is an ambient isotopy taking $S$ to a minimal surface.
We therefore investigate the coarse geometry of minimal
surfaces in $M_i$.

Set $\epsilon/2$ to be the injectivity radius of $M$.
Let $f(\kappa, \epsilon)$ be the function from Theorem 8.
Let $\tilde K$ be the lift of the 2-complex $K$
to the universal cover of $M$.
Let $k_5$ be the maximum number of complementary regions of
$\tilde K$ that lie within a distance $f(\kappa, \epsilon)$ 
of any point, and let $k_6 = f(\kappa, \epsilon) k_5$.

\noindent {\bf Lemma 13.} {\sl Let $S$ be a minimal
surface in $M_i$. Then $S$ intersects at most $k_6 |\chi(S)|$
complementary regions of $K_i$. Hence, running through
any such region, there are at most $k_6 |\chi(S)|$ translates
of $S$ under the covering group action of $G_i$.}

\noindent {\sl Proof.} By Theorem 8, the number of balls
of radius $f(\kappa, \epsilon)$
required to cover $S$ is at most $f(\kappa, \epsilon) |\chi(S)|$.
The centre of each of these balls has at most $k_5$
complementary regions of $K_i$ within a distance $f(\kappa, \epsilon)$. 
So, $S$ intersects at most $k_6 |\chi(S)|$ complementary regions of $K_i$. 
Each such region corresponds to an element of $G_i$.
To prove the second half of the lemma, we may concentrate
on the region corresponding to the identity. Then a translate
$gS$ runs through here, for some $g$ in $G_i$, if and
only if $S$ runs through the region corresponding to 
$g^{-1}$. Thus, there can be at most $k_6 |\chi(S)|$
such $g$. $\square$

\noindent {\sl Proof of Theorem 1(2).} 
Let $F$ be a strongly irreducible Heegaard surface
in $M_i$ with $|\chi(F)| = \chi_-^{sh}(M_i)$.
By Theorem 6, there
is an ambient isotopy taking it either to a minimal
surface or to the double cover of a minimal non-orientable
surface, with a small tube attached. So, 
by Lemma 13, $F$ intersects at most
$k_6 \chi_-^{sh}(M_i)$ complementary regions of $K_i$.
Hence, by Lemma 13, the number of copies of $S$ that $F$
intersects is at most $(k_6 \chi_-^{sh}(M_i))(k_6 k_4 \chi_-^h(M_i))$.
This is less than $d_i$ if $\chi_-^{sh}(M_i)/ \sqrt{d_i}$
is sufficiently small.
So there is a translate of $S$ which misses $F$. It then lies
in a complementary handlebody of $F$. But this is
impossible, since $S$ is incompressible. 
So, $\chi_-^{sh}(M_i)/ \sqrt{d_i}$ is bounded away from zero.
$\square$

\noindent {\sl Proof of Theorem 1(3).}
For ease of notation, let $x = \chi_-^h(M_i)$.
Let $\{ F_1, \dots, F_n \}$ be a
generalised Heegaard splitting for $M_i$, 
satisfying (1) - (5) of Theorem 3, obtained from 
a minimal genus Heegaard splitting. 
Replace any components of $F_1 \cup \dots \cup F_n$
that are parallel by a single component, and let
$\overline F$ be the resulting surface. Isotope
$\overline F$ so that each component is as in
Theorems 5 or 6. Corollary 4 states that
$|\chi(\overline F)| < x^2$. By Lemma 13, the number of complementary
regions of $K_i$ that can intersect $\overline F$
is at most $k_6 x^2$. Let $D$ be the corresponding
subset of $G_i$. 

Similarly, let $C$ be the subset of $G_i$ that
corresponds to those complementary regions of $K_i$
which intersect $S$. By Lemma 13, $|C| \leq k_6 |\chi(S)| \leq k_6 k_4 x$.

We claim that, when $i$ is sufficiently large, there 
are at least $9x/2$ disjoint translates
of $S$ under $G_i$ that are also disjoint from $\overline F$.
Let $m = 9x/2$. If the claim is not true, then
for any $m$-tuple $(g_1 S, \dots, g_m S)$
of copies of $S$ (where $g_j \in G_i$ for each $j$),
either at least two intersect or one copy intersects $\overline F$.
In the former case, $g_j c_1 = g_k c_2$, for some $c_1$
and $c_2$ in $C$, and for $1 \leq j < k \leq m$.
Hence, $g_k^{-1} g_j \in C C^{-1}$. In the latter
case, $g_j c_1 = d$ for some $c_1$ in $C$ and $d$ in $D$,
and so $g_j \in D C^{-1}$.
Thus, the sets $q_{jk}^{-1}(C C^{-1})$ and $p_j^{-1}(D C^{-1})$
cover $(G_i)^m$, where $q_{jk}$ and $p_j$ are the maps
$$\eqalign{
q_{jk} \colon (G_i)^m &\rightarrow G_i \cr
(g_1, \dots, g_m) &\mapsto g_k^{-1} g_j}$$
$$\eqalign{
p_j \colon (G_i)^m &\rightarrow G_i \cr
(g_1, \dots, g_m) &\mapsto g_j,}$$
for $1 \leq j < k \leq m$. 
The former sets $q_{jk}^{-1}(C C^{-1})$ 
each have size $|G_i|^{m-1} |C C^{-1}|$,
and the latter sets $p_j^{-1}(D C^{-1})$ have size $|G_i|^{m-1} |D C^{-1}|$.
So,
$$|G_i|^m \leq \left( {m \atop 2} \right ) |G_i|^{m-1} |C|^2 +
m |G_i|^{m-1} |C| |D|.$$
This implies that 
$$d_i = |G_i| \leq \left( {m \atop 2} \right ) (k_6 k_4 x)^2 + 
m (k_6 k_4 x)(k_6 x^2).$$
The right-hand side has order $x^4$ as $i \rightarrow \infty$.
However, $x/\root 4 \of {d_i} \rightarrow 0$, which is
a contradiction, proving the claim.

Consider these $9x/2$ copies of $S$.
Each lies in the complement of $\overline F$, which is
a collection of compression bodies. Since $S$ is incompressible and
connected, each copy of $S$ must be parallel to a component of $\overline F$. 
By Corollary 4(2), $\overline F$ has at most $3x/2$ components.
So, at least 3 copies of $S$ are parallel, and at least 2 of these
are coherently oriented. The proof is now completed
by the following lemma. $\square$

\noindent {\bf Lemma 14.} {\sl Let $S$ be a connected, 
embedded, oriented, incompressible, 
non-separating surface in a closed orientable 
3-manifold $M_i$. Suppose that the image of $S$ under some finite 
order orientation-preserving homeomorphism $h$ of $M_i$ is disjoint
from $S$, parallel to it and coherently oriented. Then
$M_i$ fibres over the circle with fibre $S$.}

\noindent {\sl Proof.} Let $Y$ be the manifold
lying between $S$ and $h(S)$. It is copy of $S \times I$,
with $S$ and $h(S)$ corresponding to $S \times \{ 0 \}$
and $S \times \{ 1 \}$. Take a countable collection
$\{ Y_n : n \in {\Bbb Z} \}$ of copies of this manifold.
Glue $S \times \{ 1 \}$ in $Y_n$ to $S \times \{ 0 \}$ in $Y_{n+1}$,
via $h^{-1}$. The resulting space $Y_\infty$
is a copy of $S \times {\Bbb R}$. Let $H$ be
the automorphism of this space taking $Y_n$ to $Y_{n+1}$
for each $n$, via the `identity'. Let $p \colon Y_0 \rightarrow Y$
be the identification homeomorphism.
Extend this to a map $p \colon Y_\infty \rightarrow M_i$
by defining $p|Y_n$ to be $h^n p H^{-n}$.

We claim that this is a covering map. It may
expressed as a composition $Y_\infty \rightarrow 
Y_\infty/ \langle H^N \rangle
\rightarrow M_i$, where $N$ is the order of $h$. The first
of these maps is obviously a covering map. The second is also,
since it is a local homeomorphism and $Y_\infty / \langle H^N \rangle$
is compact. Hence, $p$ is a covering map.

By construction, $h^n(S)$ lifts homeomorphically to 
$Y_{n-1} \cap Y_{n}$, for each $n$. Hence, the inverse image 
of $S$ in $Y_\infty$ includes all translates of $Y_{-1} \cap Y_0$ 
under $\langle H^N \rangle$. These translates divide
$Y_\infty$ into copies of $S \times I$. Since $p^{-1}(S)$
is incompressible and any closed embedded incompressible surface in 
$S \times I$ is horizontal, we deduce that $p^{-1}(S)$ divides $Y_\infty$ 
into a collection of copies of $S \times I$.
The restriction of $p$ to one of these components $Z$
is a covering map to a component of $M_i - S$.
But $M_i - S$ is connected, as $S$ is connected
and non-separating. So, $p$ maps $Z$ surjectively
onto $M_i - S$. By examining this map near $S$, we see that it
is degree one and hence a homeomorphism. Therefore, $M_i$ is
obtained from a copy of $S \times I$ by gluing
its boundary components homeomorphically. So, $M_i$
fibres over the circle with fibre $S$. $\square$

\vskip 12pt
\centerline {\caps References}
\vskip 6pt

\item{1.} {\caps R. Brooks}, {\sl The spectral geometry of a 
tower of coverings}, J. Differential Geom. {\bf 23} (1986) 97--107.

\item{2.} {\caps M. Freedman, J. Hass and P. Scott},
{\sl Least area incompressible surfaces in $3$-manifolds,}
Invent. Math. {\bf 71} (1983) 609--642.

\item{3.} {\caps M. Lackenby}, {\sl Heegaard splittings,
the virtually Haken conjecture and Property ($\tau$)},
Preprint.

\item{4.} {\caps A. Lubotzky}, {\sl Eigenvalues of the 
Laplacian, the first Betti number and the congruence subgroup
problem,} Ann. Math. (2) {\bf 144} (1996) 441--452.

\item{5.} {\caps J. Pitts and J. H. Rubinstein,} {\sl Existence 
of minimal surfaces of bounded topological type in three-manifolds.}
Miniconference on geometry and partial differential equations 
(Canberra, 1985), 163--176.

\item{6.} {\caps M. Scharlemann,} {\sl Heegaard splittings},
Handbook of Geometric Topology (Elsevier, 2002), 921--953.

\item{7.} {\caps M. Scharlemann and A. Thompson}, {\sl
Thin position for $3$-manifolds,} Geometric topology (Haifa,
1992), 231--238, Contemp. Math., 164.

\item{8.} {\caps B. White}, {\sl The space of minimal submanifolds
for varying Riemannian metrics}, Indiana Math. Journal {\bf 40} 
(1991) 161--200.

\vskip 12pt

\+ Mathematical Institute, Oxford University, \cr
\+ 24-29 St Giles', Oxford OX1 3LB, UK. \cr

\end